# Mathematical aspects of Multiwell Deconvolution and its relation to Capacitance Resistance Model

**Artur Aslanyan** (*Nafta College*)

E-mail: ama@naftacollege.com

## Abstract

The paper provides introduction into the mathematical aspects of Multiwell Deconvolution (MDCV) and Capacitance Resistance Model (CRM) and connection between them. Both methods are trying to train a model over the long-term history of surface flowrates and bottomhole pressure readings and then predict bottomhole and formation pressure in response to a given production/injection flowrate scenario (called "rate control simulation") or alternatively may predict flowrate and formation pressure in response to a given bottomhole pressure scenario (called "pressure control simulation"). It has been shown that CRM can be viewed as a partial case of MDCV with a specific type of a drawdown and cross-well pressure transient responses which is not always met in practice. The paper also explains limitations which are common for both methods and specify additional limitations of CRM which MDCV can handle.

# Contents

# 1. Introduction

The idea of building a correlation between flowrates and bottomhole pressure for multi-well systems is being perpetually circulated among industry specialists.

A simple and straightforward approach is suggested by the statistical regression models where one can build linear or non-linear regression (including artificial neural network) and may pass a successful blind-testing for predictions on cross-validating dataset.

Unfortunately, the statistical regression models have two major problems:

- Non-physical correlations. They often pick non-physical correlation between the wells which happens as a result of coincidence in synchronous flowrate-pressure variations. The attempt to analyze this correlation may lead petroleum engineer to a wrong understanding of well and reservoir performance.

- Non-physical predictions. When successfully trained and validated, the statistical models still may not provide a reasonable predictions outside the historical range of flowrate-pressure variations in training datasets. In other words, the validation success may not be the ultimate qualifier of the model physicality since the validation dataset is almost never representative of a broad range of possible production scenarios.

Unfortunately, there is no tool to assess if the above issues are present in the regression model until the moment when the new field data arrive. This often results in misleading predictions and wrong reservoir understanding and eventually demotivates end users.

These issues can be partially addressed within another class of correlation models which impose some limitations to flowrate-pressure behavior based on physical principles of pressure build-up and cross-well interference. The main advantage of physical models is that they can be trained over a rather limited historical dataset and yet provide a reliable prediction in a wide range of production scenarios and properly qualify reservoir performance. The most popular physical flowrate-pressure correlation models are Mutliwell Deconvolution model (**MDCV**) and Capacitance Resistance Model (**CRM**).

Although none of the two is trying to solve diffusion equation for a multi-well system but the concepts of the both models narrow down the flowrate-pressure correlation possibilities to honour some basic physical requirements and eventually will:

- minimise the risk of getting non-physical correlation

- expand the reliable model prediction range of flowrate-pressure variations far beyond the range of the historical dataset used for training the model

The main conceptual difference between **MDCV** and **CRM** is the degree of the physical conditioning.

The **CRM** is fully physical but narrows down the possible correlations to the very specific subset of diffusion equation solutions and as a consequence cannot cover numerous practical cases which fall outside this subset (see a comprehensive list of these conditions at the end of this paper).

On a contrary, the **MDCV** is based on physical constrains which cover a slightly wider range of correlations comparing to those imposed by diffusion equation and as such covers all physical possibilities of pressure-rate correlations but theoretically may possibly end up in a non-diffusion correlation.

The commercial workflow based on **MDCV** will always include the analytical or numerical diffusion modelling of the **MDCV** output (deconvolved pressure) which eventually provides the fully physical pressure-rate correlation.

The trick is that analytical (or numerical) diffusion modelling of the **MDCV** output (called Unit-rate Transient Responses, see below) is a very simple, fast and straightforward procedure, comparing to the attempt of tuning up the diffusion model to match the original historical dataset. One may see **MDCV** as a mathematical pre-processing to simplify the diffusion modelling and as such the **MDCV**-based workflow always ends up in a fully-physical correlation.

One can get some high-level readings on **MDCV** in [7, 8, 12, 14, 19, 20, 22, 23].
The practical applications of **MDCV** can be found in [1 - 5], where it is listed as the mathematical engine of the Multi-well Retrospective Testing (**MRT**) service. The Capacitance Resistance Model (**CRM**) is a popular Production Analysis tool and has been around for few decades already [11, 16, 21].

Despite of its mathematical simplicity, long presence on the market and numerous practical trials the **CRM** model did not become the common industry tool due to a frequent failure to deliver a reliable forecast. In author's practise, based on response from the numerous petroleum companies, the majority of **CRM** failures were not properly explained by users, or were not explained at all. This tendency has created a clear call for deeper study and understanding of when and how **CRM** should actually be applied.

As it turned out of the author's investigation, many **CRM** failure cases were not related to the capability of **CRM** engine but rather to the misunderstanding and mismanagement of the input data and violation of cross-validation procedure. As for the cases when the data and workflow were fair and **CRM** still failed then one of the main reasons, as it will come clear below, was a violation of the fundamental conditions at which **CRM** model has been derived. Some of them are cardinal and none of the flowrate-pressure models can overcome them (like poor flowrate variations). While the others (and covering the majority of cases) can be handled by the more flexible flowrate-pressure correlation models, like Multiwell Deconvolution (**MDCV**).

Just as with statistical regression models, the key metric of the **MDCV** and **CRM** model success is the ability to pass a blind cross-validation. In order to do this, the original historical dataset is split into two subsets: training dataset and validation dataset.

The training dataset is used to train the model while the validation dataset is used to assess the predictive quality of the trained mode, normally presented as the cross-plot of prediction vs actual data and quantified by Root Mean Square Deviation (RMSD) and Coefficient of Determination (R2).

The selection of the validation dataset is done by the analyst and constitutes the important input in this analysis. The procedure of selecting the validation dataset is quite straightforward but in some practical cases the production history is uneventful, and this will make analyst struggle to select a truly representative validation dataset.

The simplest way to understand if this validation dataset is good for the job is to build an express numerical pressure diffusion 2D+ model ("plus" meaning multilayer case) and see if a given validation dataset can discriminate between accurate and inaccurate solutions (with accuracy being defined by the analyst depending on the context of his study).
Still this approach is not the ultimate judge as it's still based on some propositions on reservoir properties which have been used in 2D+ sector model but it's better than measure success with numbers borrowed from other studies. It should be mentioned that in academic studies the 2D+ sector models work really well to analyse the performance of different flowrate-pressure models and summarise their scope of applicability.

The present study provides the derivation of **MDCV** and **CRM** driving equations, the mathematical and physical link between model parameters and high-level scope of applicability of both models.

## 2. Derivation of MDCV equations

In many practical cases the multiphase fluid reservoir flow, created by production from or injection to a group of $N$ wells, will honour a linear pressure diffusion equation for slightly compressible fluid:

$$\beta \frac{\partial p}{\partial t} = \nabla \left( \widehat{M}(\nabla p - \rho_0 \cdot \mathbf{g}) \right) + \sum_n q_n(t) \delta(\mathbf{r} - \mathbf{r}_n) \tag{2.1}$$

with initial conditions:

$$p(t = 0, \mathbf{r}) = p_0(\mathbf{r}) \tag{2.2}$$

$$\widehat{M}(\nabla p_0(\mathbf{r}) - \rho_0 \cdot \mathbf{g})\big|_{t=0} = 0 \tag{2.3}$$

and general form of reservoir boundary $\mathbf{\Gamma}$ condition:

$$\left[ a \cdot (p(t, \mathbf{r}) - p_0) + \epsilon \cdot \mathbf{n} \cdot \widehat{M}(\nabla p(t, \mathbf{r}) - \rho_0 \cdot \mathbf{g}) \right]_{r \in \Gamma} = 0, \quad a^2 + \epsilon^2 > 0 \tag{2.4}$$

where:

$p(t, \mathbf{r}) = \frac{1}{3}(p_w + p_o + p_g)$ is 3-phase-average reservoir fluid pressure,

$p_0(\mathbf{r}) = p(t = 0, \mathbf{r})$ is initial pressure distribution across the reservoir,

$\phi(\mathbf{r}) = \phi^*(\mathbf{r}) \cdot \exp[c_r(p) \cdot (p - p^*)]|_{p=p_0(\mathbf{r})}$ is reservoir porosity at initial reservoir pressure,

$\phi^*(\mathbf{r})$ is reservoir porosity at the reference reservoir pressure $p^*$ (usually selected as a hydrostatic pressure at datum),

$\{s_w(\mathbf{r}), s_o(\mathbf{r}), s_g(\mathbf{r})\}$ are water, oil and gas reservoir saturations, subject to $s_w + s_o + s_g = 1$ and not changing in time,

$\rho_0(\mathbf{r}) = \rho_w(\mathbf{r}) \cdot s_w(\mathbf{r}) + \rho_o(\mathbf{r}) \cdot s_o(\mathbf{r}) + \rho_g(\mathbf{r}) \cdot s_g(\mathbf{r})$ is total fluid density,

$\{\rho_w(\mathbf{r}) = \rho_w(p), \rho_o(\mathbf{r}) = \rho_o(p), \rho_g(\mathbf{r}) = \rho_g(p)\}|_{p=p_0(\mathbf{r})}$ are water, oil and gas densities calculated at initial pressure conditions,

$c_t(\mathbf{r}) = c_r(\mathbf{r}) + c_w(\mathbf{r}) \cdot s_w(\mathbf{r}) + c_o(\mathbf{r}) \cdot s_o + c_g(\mathbf{r}) \cdot s_g(\mathbf{r})$ is total reservoir compressibility,

$c_r(\mathbf{r}) = \frac{1}{\phi} \frac{\partial \phi}{\partial \tilde{p}}\big|_{p=p_0(\mathbf{r})}$ is pore compressibility,

$\{c_w(\mathbf{r}) = c_w(p), c_o(\mathbf{r}) = c_o(p), c_g(\mathbf{r}) = c_g(p)\}|_{p=p_0(\mathbf{r})}$ are water, oil and gas compressibility calculated at initial pressure conditions,

$\beta(\mathbf{r}) = c_t(\mathbf{r}) \cdot \phi(\mathbf{r})$ is reservoir storativity,

$\hat{k}_a(\mathbf{r})$ is non-degenerate symmetric tensor of absolute reservoir permeability,

$\{\mu_w(\mathbf{r}) = \mu_w(p), \mu_o(\mathbf{r}) = \mu_o(p), \mu_g(\mathbf{r}) = \mu_g(p)\}|_{p=p_0(\mathbf{r})}$ are water, oil and gas dynamic viscosities calculated at initial pressure conditions,

$\{k_{rw}, k_{ro}, k_{rg}\}$ are water, oil and gas relative permeabilities as functions of reservoir saturation $\{s_w, s_o, s_g\}$,

$\widehat{M} = \left[\frac{k_{rw}}{\mu_w} + \frac{k_{ro}}{\mu_o} + \frac{k_{rg}}{\mu_g}\right] \cdot \hat{k}_a$ is a tensor of reservoir fluid mobility,

$\mathbf{g} = g \cdot \mathbf{n}_g$ is standard gravity vector,

$\mathbf{n}_g$ is vector pointing downward to the Earth's gravity center,

$g = 9.81$ m²/s is Earth's standard gravity,

$a(\mathbf{r}) \geq 0$ is some non-negative function,

$\epsilon(\mathbf{r}) \in \{0,1\}$ is function which binary values (either 0 or 1),

$\mathbf{n}$ is external normal to the reservoir boundary $\Gamma$,

$\nabla f = \frac{\partial f}{\partial x^i}$ is pressure gradient of a scalar function $f(t, \mathbf{r})$ along the position vector $\mathbf{r} = \{x^1, x^2, x^3\} \in \mathrm{R}^3$,

$\nabla \mathbf{b} = \frac{\partial b_1}{\partial x^1} + \frac{\partial b_2}{\partial x^2} + \frac{\partial b_3}{\partial x^3}$ is divergence of vector field $\mathbf{b}$,

$\widehat{M} \mathbf{b}$ is a vector defined by linear operator $\widehat{M}$ acting on a vector field $\mathbf{b}$ with the following result:

$(\widehat{M} \mathbf{b})_i = \sum_j M_{ij} \cdot b_j, \quad i, j \in \{1,2,3\}$,

$\{q_n(t)\}_{n \in \{1..N\}}$ are flowrate histories with a natural condition:

$$q_n(t \leq 0) = 0, \forall n \in \{1..N\} \qquad (2.5)$$

The production rates in the above convention are positive: $q_{prod}(t) = q^{\uparrow}(t) > 0$ and injection rates are negative: $q_{inj}(t) = -q^{\downarrow}(t) < 0$.

The diffusion model (2.1) – (2.4) assumes that:

- reservoir saturation does not vary in time:
  $\{s_w(t) = s_w = \text{const}, s_o(t) = s_o = \text{const}, s_g(t) = s_g = \text{const}\}$
- reservoir pressure variation is low enough to assume that absolute permeability $\hat{k}_a(\mathbf{r})$, porosity $\phi$, total compressibility $c_t$ and fluid viscosities $\{\mu_w, \mu_o, \mu_g\}$ do not depend on pressure.

As result of the above assumptions the storativity $\beta$ and mobility $\widehat{M}$ will be not dependent on time and pressure and the boundary value problem (2.1)–(2.4) is represented by a set of linear equations and linear initial and boundary conditions.

The condition (2.3) is also called a "stationary start" and means that initially there are no reservoir flows and it will stay the same unless production/injection starts. The reservoir boundary condition (2.4) is very general and covers the most (if not all) practical cases, including the "constant pressure support" flow ($a = 1, \epsilon = 0$), fully restricted "no flow" condition ($a = 0, \epsilon = 1$) and partial pressure support (as from aquifer or gas cap) ($a > 0, \epsilon = 1$).

The boundary value problem (2.1)–(2.4) represents a Cauchy problem with the unique solution. It should be noted that (2.1)–(2.4) assumes that coefficients are constant in time which may not hold true if reservoir properties (like for example saturation) are changing over time. This is usually a slow process (taking months and years) and can be manually split into time-segments where (2.1)–(2.4) is a relevant approximation. In general case of non-linear multiphase flow the phase-average pressure may not honour equations (2.1)–(2.4) even at short times.

There are many tricks to tackle this problem.

The non-linear effects related to fluid compressibility (including multi-phase fluids and dry gas) can be substantially compensated by use of pseudo-variables (pseudo-pressure and pseudo-time).

The non-darcy flow effects (usually present in high rate gas wells) are mostly affecting the well vicinity and as such can be approximated by rate-dependant skin-factor within time interval where rate does not have substantial step changes. Discussing these methods falls outside the scope of this paper but overall one can conclude that in the majority of practical cases the pressure diffusion may be constructed by means of solutions of equations (2.1)–(2.4) with sufficient (for engineering purposes) degree of accuracy. The derivation of (2.1)–(2.4) from general 3-phase Black Oil reservoir flow equations is quite straightforward but discussion of its applicability and limitations extends beyond the purpose of this paper. Most probably, the first step in this derivation was done by Perrine [17] followed by a more thorough work of Martin [15].

The next step towards the derivation of **MDCV** is to remove the gravity term by introducing a new pressure variable $\tilde{p}(t, \mathbf{r})$:

$$p(t,r) \to \tilde{p}(t, \mathbf{r}) = p_0(\mathbf{r}) - p(t, \mathbf{r}) \tag{2.6}$$

which transforms the boundary value problem (2.1)–(2.4) into:

$$\beta \frac{\partial \tilde{p}}{\partial t} = \nabla(\widehat{M} \nabla \tilde{p}) + \sum_n q_n(t) \delta(\mathbf{r} - \mathbf{r_n}) \tag{2.7}$$

with initial conditions:

$$\tilde{p}(t = 0, \mathbf{r}) = 0 \tag{2.8}$$

and boundary conditions:

$$[a \cdot \tilde{p}(t, \mathbf{r}) + \epsilon \cdot \mathbf{n} \cdot \widehat{M}(\nabla \tilde{p}(t, \mathbf{r})]_{r \in \Gamma} = 0, \qquad a^2 + \epsilon^2 > 0 \tag{2.9}$$

The new boundary value problem (2.7)–(2.9) will constitute the basis for further calculations.

The derivation of convolution principle starts with considering a following boundary value problem:

$$\beta \frac{\partial p_{u,m}}{\partial t} = \nabla(\widehat{M} \nabla p_{u,m}) + \theta(t) \cdot \delta(\mathbf{r} - \mathbf{r}_m), \quad p_{u,m}(t = 0, \mathbf{r}) = 0,$$

$$[a \cdot p_{u,m}(t, \mathbf{r}) + \epsilon \cdot \mathbf{n} \cdot \widehat{M} \nabla p_m(t, \mathbf{r})]_{r \in \Gamma} = 0 \tag{2.10}$$

which represents a reservoir pressure response $p_{u,m}(t, \mathbf{r})$ to the unit-rate production of $m$-th well $q_{u,m}(t) = \theta(t)$ (all other wells being shut-down) with zero initial pressure $p_{u,m}(t = 0, \mathbf{r}) = 0$, where $\theta(t)$ is Heaviside step function:

$$\theta(t) = \begin{cases} 0, & if\ t \leq 0 \\ 1, & if\ t > 0 \end{cases} \tag{2.11}$$

Consider a convolution of unit-rate reservoir pressure responses $p_{u,m}(t, \mathbf{r})$ with arbitrary flowrate variation history of $m$-th well $q_m(t)$:

$$p_m(t,\mathbf{r}) = \int_0^t p_{u,m}(t-\tau,\mathbf{r})dq_m(\tau) = \int_0^t p_{u,m}(t-\tau,\mathbf{r})\dot{q}_m(\tau)d\tau \qquad (2.12)$$

where dot $\dot{g}(\xi) = \frac{dg}{d\xi}$ means derivative with respect to the whole argument.

The time derivative of the (2.12) is going to be:

$$\frac{\partial p_m(t,\mathbf{r})}{\partial t} = \frac{\partial}{\partial t}\int_0^t p_{u,m}(t-\tau,\mathbf{r})\dot{q}_m(\tau)d\tau = \int_0^t \frac{\partial p_{u,m}(t-\tau,\mathbf{r})}{\partial t}\dot{q}_m(\tau)d\tau + \left[p_{u,m}(t-\tau,\mathbf{r})\dot{q}_m(\tau)\right]\Big|_{\tau=t} =$$
$$\int_0^t \frac{\partial p_{u,m}(t-\tau,\mathbf{r})}{\partial t}\dot{q}_m(\tau)d\tau + p_{u,m}(0,\mathbf{r})\dot{q}_m(\tau) = \int_0^t \frac{\partial p_{u,m}(t-\tau,\mathbf{r})}{\partial t}\dot{q}_m(\tau)d\tau \qquad (2.13)$$

where $p_{u,m}(0,\mathbf{r}) = 0$ is the initial condition for the unit-rate reservoir pressure response (2.12).

By convolving the equations (2.10) with the well rate variation $q_m(t)$ and taking into account the equation (2.13) one arrives to conclusion that (2.12) satisfies the following equations:

$$\beta \frac{\partial p_m}{\partial t} = \nabla(\widehat{M}\nabla p_m) + q_m(t)\cdot\delta(\mathbf{r}-\mathbf{r}_m), \quad p_m(t=0,\mathbf{r}) = 0,$$

$$[a\cdot p_m(t,\mathbf{r}) + \epsilon\cdot\mathbf{n}\cdot\widehat{M}\nabla p_m(t,\mathbf{r})]_{r\in\Gamma} = 0 \qquad (2.14)$$

providing that storativity $\beta$ and mobility $\widehat{M}$ are not varying in time during the modelling period.

The solution of the boundary value problem (2.14) represents a reservoir pressure response to the $m$-th well production/injection with the arbitrary flowrate history $q_m(t)$.

Now consider a linear combination of pressure responses $\{p_m(t,\mathbf{r})\}_{m\in\{1..N\}}$ to production/injection in all wells:

$$\tilde{p}(t,\mathbf{r}) = \sum_m p_m(t,\mathbf{r}) = \sum_m \int_0^t p_{u,m}(t-\tau,\mathbf{r})dq_m(\tau) \qquad (2.15)$$

By direct substitution one can see that (2.15) satisfies all the equations (2.7)–(2.9).

Consequently, a difference between initial pressure distribution $p_0(\mathbf{r})$ and $\tilde{p}(t,\mathbf{r})$:

$$p(t,\mathbf{r}) = p_0(\mathbf{r}) - \tilde{p}(t,\mathbf{r}) = p_0(\mathbf{r}) - \sum_m p_m(t,\mathbf{r}) = p_0(\mathbf{r}) - \sum_m \int_0^t p_{u,m}(t-\tau,\mathbf{r})dq_m(\tau) \qquad (2.16)$$

is going to satisfy all equations (2.1)–(2.4) and hence represents the unique solution to the original boundary value problem for $N$ wells with the arbitrary flowrate histories $\{q_m(t)\}_{m\in\{1..N\}}$.

Consider solution of (2.16) in the location $\mathbf{r}_n$ of the $n$-th well: $p_n(t) = p(t,\mathbf{r}_n)$ which represents a bottomhole pressure response in $n$-th well to the production/injection activity in the same and offset wells.

Let's define Drawdown Transient Response = **DTR** (sometimes called **self-response**) as $p_{u,nn}(\tau) = p_{u,n}(\tau,\mathbf{r}_n)$ representing a pressure response in location $\mathbf{r}_n$ of the $n$-th well to the unit-rate production in the same well $q_{u,n}(\tau) = \{0, \tau \leq 0 | 1, \tau > 0\}$.

Let's define Cross-well Transient Response = **CTR** (sometimes called **interference-response**) as $p_{u,nm}(\tau) = p_{u,m}(\tau,\mathbf{r}_n)$ representing a pressure response in location $\mathbf{r}_n$ of the $n$-th well to the unit-rate production in the offset $m$-th ($m\neq n$) well $q_{u,m}(\tau) = \{0, \tau \leq 0 | 1, \tau > 0\}$.

The **DTR** and **CTR** are also called Unit-Rate Transient Response (**UTR**) and play a key role in deconvolution theory.

The equation (2.16) can be now rewritten in terms of the bottomhole pressure of the $n$-th well:

$$p_n(t) = p_{0n} - \sum_m \int_0^t p_{u,nm}(t-\tau)dq_m(\tau) =$$

$$= p_{0n} - \int_0^t p_{u,nn}(t-\tau)dq_n(\tau) - \sum_{m \neq n} \int_0^t p_{u,nm}(t-\tau)dq_m(\tau) \quad (2.17)$$

where $p_{0n} = p_n(t=0)$ as $q_m(t) = 0, \forall t \in (-\infty, 0], \forall m \in \{1..N\}$ (see (2.5)).

Equation (2.17) is called a **Multiwell Convolution Equation** for bottomhole pressure and means that if all **DTR/CTR** are available then a pressure response in each well $p_n(t)$ can be presented as a linear sum of initial formation pressure $p_{0n}$, a convolution of **DTR** $p_{u,nn}(t)$ with $n$-th well flowrate history $q_n(t)$ and sum of all convolutions of **CTRs** $p_{u,nm}(t)$ with offset wells flowrate variation histories $\{q_m(\tau)\}_{m \in \{1..N | m \neq n\}}$.

We shall further derive the alternative form of the **Multiwell Convolution Equation** which is very useful in many calculations.

Let's define the derivative with respect to the whole argument as $\dot{g}(\xi) = \frac{dg}{d\xi}$ so that $dq_n(\tau) = \dot{q}_n \cdot d\tau$ and (2.17) will take a form:

$$p_n(t) = p_{n0} - \int_0^t p_{u,nn}(t-\tau)\dot{q}_n(\tau)d\tau - \sum_{m \neq n} \int_0^t p_{u,m}(t-\tau)\dot{q}_m(\tau)d\tau \quad (2.18)$$

Note the property of convolution:

$$-\int_0^t \dot{a}(t-\tau)b(\tau)d\tau + \int_0^t a(t-\tau)\dot{b}(\tau)d\tau = [a(t-\tau)b(\tau)]|_0^t = a(0)b(t) - a(t)b(0) \quad (2.19)$$

Applying this to the first integral in (2.18):

$$\int_0^t p_{u,nn}(t-\tau)\dot{q}_n(\tau)d\tau = p_{u,nn}(0)q_n(t) - p_{u,nn}(t)q_n(0) + \int_0^t \dot{p}_{u,nn}(t-\tau)q_n(\tau)d\tau \quad (2.20)$$

and since $p_{u,nn}(0) = 0$ and $q_n(0) = 0$ then:

$$\int_0^t p_{u,nn}(t-\tau)\dot{q}_n(\tau)d\tau = \int_0^t \dot{p}_{u,nn}(t-\tau)q_n(\tau)d\tau \quad (2.21)$$

so that (2.18) will take a form:

$$p_n(t) = p_{n0} - \int_0^t \dot{p}_{u,nn}(t-\tau)q_n(\tau)d\tau - \sum_{m \neq n} \int_0^t p_{u,nm}(t-\tau)\dot{q}_m(\tau)d\tau \quad (2.22)$$

The similar exercise with the last term:

$$\int_0^t p_{u,nm}(t-\tau)\dot{q}_m(\tau)d\tau = p_{u,nm}(0)q_m(t) - p_{u,nm}(t)q_m(0) + \int_0^t \dot{p}_{u,nm}(t-\tau)q_m(\tau)d\tau =$$

$$\int_0^t \dot{p}_{u,nm}(t-\tau)q_m(\tau)d\tau \quad (2.23)$$

so that:

$$p_n(t) = p_{n0} - \int_0^t \dot{p}_{u,nn}(t-\tau)q_n(\tau)d\tau - \sum_{m \neq n} \int_0^t \dot{p}_{u,nm}(t-\tau)q_m(\tau)d\tau \quad (2.24)$$

The above equation (24) is the equivalent form of the original convolution equation (2.17).

The main goal of **Multiwell Deconvolution** (**MDCV**) is to solve equation (2.17) or (2.24) by finding the initial formation pressure $p_{n0}$ and a set of UTRs ($\{p_{u,nm}(t)\}_{n,m \in \{1..N\}}$) based on the available field data on flowrate histories $\{q_n(\tau)\}_{n \in \{1..N\}}$ and bottomhole pressure histories $\{p_n(\tau)\}_{n=1..N}$. Sometimes the initial formation pressure $p_{n0}$ is known and this simplifies the **MDCV** procedure.

The main assumption of **MDCV** is that diffusion process honors linear equations (2.1)–(2.4) and as a consequence the UTRs stay constant during the modelling period. There are some specialised **MDCV** techniques which use linear deconvolution for non-linear pressure diffusion by means of pseudo-potentials and varying skin. There is also a possibility to generalise the (2.17) for the non-stationary start with $p_{0n}(t)$ being dependent on time but quickly decaying (as a result of the unknown pre-start production activity). These above modifications extend the applicability of **MDCV** engine to many practical cases which do not directly honour the equations (2.1)–(2.4). Discussing these methods falls outside the scope of this paper and does not change the main conclusions below including those related to the scope of **MDCV** applicability and comparison against **CRM**.

The general idea behind **MDCV** engine is to create an optimisation loop over initial formation pressure $p_{n0}$ (in case it is not known), **DTR/CTR** $\{p_{u,nm}\}_{n,m \in \{1..N\}}$ and flowrate corrections $\{\delta q_n(t) = q_n(t) - \tilde{q}_n(t)\}_{n \in \{1..N\}}$ (in case they are not accurately allocated) to minimize the objective function:

$$E[p_0, \{p_{u,nm}\}, q_n] = \sum_n \sum_{\alpha_p} w_{p,n\alpha_p}(p_n(t_{\alpha_p}) - \tilde{p}_n(t_{\alpha_p}))^2 + \sum_n \sum_{\beta_q} w_{q,n\beta_q}(q_n(t_{\beta_q}) - \tilde{q}_n(t_{\beta_q}))^2 \to min \quad (2.25)$$

where $\{p_n(t), q_n(t)\}$ is **MDCV** model realisation, $\{\tilde{p}_n(t_{\alpha_p}), \tilde{q}_n(t_{\beta_q})\}$ is input data set on bottomhole pressure and flowrate histories for all tested wells recorded at discrete time moments $\{t_\alpha\}$ and $\{w_{p,n\alpha}, w_{q,n\beta}\}$ is a set of trust weights for each data points.

From mathematical point of view, solving equation (2.17) in general form is a typical ill-posed problem, meaning that it may have infinite number of different solutions honoring the same flowrate-pressure dataset.

Developing a proper regularization procedure is the key to success. The optimisation iterations are usually performed in a smaller functional space to improve the convergence, filter out non-physical solutions and reduce uncertainty in the final solution and at the same time avoiding substantial limitations to the possible behavior of transient responses so that all practically important solutions of boundary value problem (2.1)–(2.4) can be picked up by the deconvolution. This can be done explicitly by expanding the variables $\{p_{u,nm}\}_{n,m=1..N}$ and $\{q_n\}_{n=1..N}$ with respect to some base functions constituting the limited search space.

Alternatively this can be done by adding additional terms into the objective function (2.25), like, for example, a **UTR** curvature, deviation between model and historical rate cumulatives, deviation between model rates and well tests. The process should be also stable to typical data contamination, which requires a specific data pre-processing (like wavelet filtering and automatic pressure-rate synchronization) and specialised optimisation algorithms (usually a hybrid of Differential Evolution [9] and customised Gauss-Newton [10] algorithms). The commercial algorithm will be most probably using a combination of all the above methods.

Finally, once deconvolution algorithm arrived in some solution it should be validated on previously selected validation data subset and in case it does not match as accurately as expected it should automatically change the objective function weights and other regularisation parameters and repeat the optimisation process again.

Once the validated solution is achieved it should undergo a sensitivity analysis to assess the uniqueness of this solution. Splitting initial dataset into a training and validation datasets is usually done by the user based on the physical meaning of the production history and the amount and the contrast of rate variation events.

Apart from that, some commercial algorithms perform inner split of the training dataset into training-validation subsets (called bootstrapping [6]) in order to auto-tune the parameters of the goal function and governing properties of optimization loop to achieve the ultimate convergence or report a failure without infinite looping.

Developing a stable deconvolution engine is a substantial mathematical challenge and requires solid efforts and extended coding to achieve practically acceptable results, so that (2.25) gives only a high level impression of the deconvolution algorithm leaving behind the massive amount of technicalities, which are falling out the scope of this paper.

## 3. Derivation of CRM equations

Conventional Capacitance Resistance Model (**CRM**) model defines the flowrate response $q_n^\uparrow(t)$ in $n$-th producer to the variation of its bottomhole pressure $p_n(t)$ and the offset injection $q_m^\downarrow(t)$:

$$q_n^\uparrow(t) + \tau_n \cdot \frac{dq_n^\uparrow}{dt} = \sum_m f_{nm} \cdot q_m^\downarrow(t) - \gamma_n \cdot \frac{dp_n}{dt} \qquad (3.1)$$

where:

$q_n^\uparrow(t)$ is production rate from $n$-th producer,

$p_n(t)$ is bottomhole pressure $n$-th producer,

index $_m$ is running through injectors only: $m \in \downarrow$.

$q_m^\downarrow(t)$ is injection rate in $m$-th injector,

$f_{nm}$ is a fraction of total injection in $m$-th injector which supports production from $n$-th producer,

$\tau_n = \frac{c_{t,n} V_{\phi,n}}{J_n}$, $\gamma_n = c_{t,n} V_{\phi,n}$ are the model parameters related to $n$-th producer and reservoir properties around it,

$c_{t,n}$ is total reservoir compressibility within the drainage volume of $n$-th producer,

$V_{\phi,n}$ is drainage pore volume of $n$-th producer,

$J_n = \frac{\gamma_n}{\tau_n}$ is productivity index of the $n$-th producer.

The model parameters are subject to the obvious physical conditions:

$$\tau_n \geq 0, \ \gamma_n \geq 0, \quad 0 \leq f_{nm} \leq 1, \quad 0 \leq \sum_n f_{nm} \leq 1 \qquad (3.2)$$

In some cases, one can impose a stronger condition on $f_{nm}$:

$$\sum_n f_{nm} = 1 \qquad (3.3)$$

which means that total injection in $m$-th injector is fully distributed between all producers and no part of it is missing outside the oil pay.

Obviously, one has to rely on surveillance data or some other source of confidence to support this proposition. In author's experience a thief injection is a very common occurrence and, in this case, the (3.3) constraint should be avoided.

**CRM** model assumes that each producer is draining its own pore volume $V_{\phi,n}$ for prolonged period of time. It may happen if a given producer drains fluid from its own reservoir compartment which is isolated from other producers. It may also happen because a given produce establishes a parity in sharing the total pore volume with other producers and $V_{\phi,n}$ is a stabilised "drainage volume share cut" of the $n$-th producer.

At the same time **CRM** model assumes that no matter of location, any injector can possibly support production in any producer, and this is defined by $f_{nm}$ constants.

The equation (3.1) can be explicitly integrated:

$$q_n^\uparrow(t) = \exp(-t/\tau_n) \cdot \left[q_n^\uparrow(0) + \tau_n^{-1} \cdot \int_0^t \exp(s/\tau_n) \left[\sum_m f_{nm} \cdot q_m^\downarrow(s) - \gamma_n \frac{dp_n}{ds}\right] ds\right] \quad (3.4)$$

and used to minimize the deviation from the historical production rates records $\{\tilde{q}^\uparrow(t_k)\}$ by searching the set of model parameters $\{\tau_n, \gamma_n, f_{nm}\}$ with following objective function:

$$E[\tau_n, \gamma_n, f_{nm}] = \sum_k [q^\uparrow(t_k) - \tilde{q}^\uparrow(t_k)]^2 \to min \quad (3.5)$$

where $\{t_k\}$ is the discrete time scale of the flowrate records.

Once model parameters are estimated then one can use equation (3.1) to reconstruct the history of formation pressure as:

$$p_{r,n}(t) = p_n(t) + J_n^{-1} \cdot q_n^\uparrow(t) \quad (3.6)$$

One can also use equation (3.1) and (3.6) to predict production performance $q_n^\uparrow(t)$ and formation pressure $p_{r,n}(t)$ of the $n$-th producer based on the future bottomhole pressure $p_n(t)$ (defined by lift settings) and offset injection $\{q_m^\downarrow(t)\}$.

The original **CRM** concept can be slightly modified to answer a different question: how does the bottomhole pressure $p_n(t)$ of the $n$-th producer respond to its flowrate history $q_n^\uparrow(t)$ with account of the offset injection $\{q_m^\downarrow(t)\}$ history, which is very close to the original formulation of **MDCV** problem.

In order to answer that question, lets integrate both sides of the original **CRM** equation (3.1):

$$\int_0^t q_n^\uparrow(s)ds + \tau_n \cdot \int_0^t \frac{dq_n^\uparrow}{ds} ds = \sum_m f_{nm} \cdot \int_0^t q_m^\downarrow(s) - \gamma_n \cdot \int_0^t \frac{dp_n}{ds} ds \quad (3.7)$$

$$Q_n^\uparrow(t) + \tau_n \cdot (q_n^\uparrow(t) - q_n^\uparrow(0)) = \sum_m f_{nm} \cdot Q_m^\downarrow(t) - \gamma_n \cdot (p_n(t) - p_0(t)) \quad (3.8)$$

and then one finally gets an explicit form for the bottom-hole pressure **CRM** prediction:

$$p_n(t) = p_n(0) - \tau_n/\gamma_n \cdot [q_n^\uparrow(t) - q_n^\uparrow(0)] - \gamma_n^{-1} \cdot Q_n^\uparrow(t) + \gamma_n^{-1} \cdot \sum_m f_{nm} \cdot Q_m^\downarrow(t) \quad (3.9)$$

where
$Q_n^\uparrow(t) = \int_0^t q_n^\uparrow(s)ds$ is cumulative production from the $n$-th producer,
$Q_m^\downarrow(t) = \int_0^t q_m^\downarrow(s)ds$ is cumulative injection to the $m$-th injector.

The equation (3.9) can be used to minimize the deviation from the historical bottomhole pressure records $\tilde{p}_n(t)$ by searching the set of model parameters $\{\tau_n, \gamma_n, f_{nm}\}$ with following objective function:

$$E[\tau_n, \gamma_n, f_{nm}] = \sum_k \sum_n [p_n(t_k) - \tilde{p}_n(t_k)]^2 \to min \qquad (3.10)$$

The equation (3.8) leads to what is called the Integrated Capacitance Resistance Model (**ICRM**) for the $n$-th producer (see [16]):

$$\tau_n \frac{dQ_n^\uparrow(t)}{dt} + Q_n^\uparrow(\tau) = \sum_{m \in \downarrow} f_{nm} Q_m^\downarrow(\tau) + \gamma_n [p_{n,0} - p_{n,wf}(t)] \qquad (3.11)$$

where:

$p_{wf}(t)$ is bottomhole pressure history of $n$-th producer,

$q_n^\uparrow(t)$ is production rate history of $n$-th producer,

$Q_n^\uparrow(t) = \int_0^t q_n^\uparrow(\tau)d\tau$ is cumulative production rate history of $n$-th producer,

$q_m^\downarrow(t)$ is injection rate history of $m$-th injector,

index $m$ is running through injectors only: $m \in \downarrow$.

$\{\tau_n, \gamma_n, f_{nm}\}$– a set of **ICRM**-model constants.

The **ICRM** -model constants are estimated through the fitting procedure with the following objective function:

$$E[\tau_n, \gamma_n, f_{nm}] = \sum_m \sum \beta_q [Q_n^\uparrow(t_{\beta_q}) - \tilde{Q}_n^\uparrow(t_{\beta_q})]^2 \to min \qquad (3.12)$$

and constraints:

$$\tau_n \geq 0, \ \gamma_n \geq 0, \quad f_{nm} \geq 0, \quad \sum_{m \in \downarrow} f_{nm} \leq 1 \qquad (3.13)$$

where

$t_{\beta_q}$ are time moments when fields record on bottomhole pressure and flowrate histories were recorded,

$\tilde{q}_n^\uparrow(t)$ is the field record of production rate history of $n$-th producer,

$\tilde{Q}_n^\uparrow(t) = \int_0^t \tilde{q}_n^\uparrow(\tau)d\tau$ is the field record of cumulative production rate history of $n$-th producer.

There are commercial **ICRM** engines which train the model over the weighted goal function combining the (3.5) and (3.12).

The other forms of **CRM** are mathematically equivalent to the **ICRM** and one can find some of them in [11] and [21].

Although mathematically equivalent they are quite different in terms of computational efficiency and the comparative analysis falls outside the scope of this paper. Further down all references to **CRM** will actually mean **ICRM** realization but all the conclusions and scope of applicability is equally applied to the other forms of **CRM**.

Let's first derive the **CRM** for one producer supported by one injector. The first assumption of **CRM** is that productivity index $J$ of producer stays constant in time:

$$J = \frac{q^\uparrow(t)}{p_r(t) - p_{wf}(t)} = \text{const} \qquad (3.14)$$

where $q^\uparrow(t)$ is production rate history, $p_{wf}(t)$ is bottomhole pressure history, $p_r(t)$ is formation pressure averaged over the drainage area around a producer.

This can be re-written as explicit formula for drainage area formation pressure:

$$p_r(t) = p_{wf}(t) + J^{-1} \cdot q^\uparrow(t) \tag{3.15}$$

The second assumption is that drainage volume of a producer is finite and constant in time:

$$V_\phi(t) = V_r \phi = \text{const} \tag{3.16}$$

The third assumption is that total reservoir compressibility stays constant in time:

$$c_t \equiv \frac{1}{V_\phi} \cdot \frac{dV_\phi}{dp} = \text{const} \tag{3.17}$$

The physical meaning of these assumptions will be discussed below.

Equation (3.17) can be rewritten as:

$$dV_\phi = c_t V_\phi dp \tag{3.18}$$

meaning that drainage volume variation $\delta V_\phi$ leads to formation pressure variation:

$$dV_\phi = c_t V_\phi \delta p = c_t V_\phi [p_0 - p_r(t)] \tag{3.19}$$

The variation in drainage volume $\delta V_\phi$ is caused by production/injection:

$$\delta V_\phi = Q^\uparrow(t) - f \cdot Q^\downarrow(t) = \int_0^t q^\uparrow(\tau) d\tau - f \cdot \int_0^t q^\downarrow(\tau) d\tau \tag{3.20}$$

where $Q^\uparrow(t)$ is cumulative production rate history, $Q^\downarrow(t)$ is cumulative injection rate history, $q^\uparrow(t)$ is instantaneous production rate history, $q^\downarrow(t)$ is instantaneous injection rate history, $f$ is a part of injection which is shared to a drainage volume of producer.

This allows rewriting (3.18) as:

$$Q^\uparrow(t) - f \cdot Q^\downarrow(t) = c_t V_\phi [p_0 - p_r(t)] \tag{3.21}$$

Substituting $p_r(t)$ from productivity equation (40) one arrives to:

$$Q^\uparrow(\tau) = f \cdot Q^\downarrow(\tau) + c_t V_\phi [p_0 - p_{wf}(t) + J^{-1} q^\uparrow(t)] \tag{3.22}$$

which leads to:

$$\tau \cdot \frac{dQ^\uparrow(t)}{dt} + Q^\uparrow(\tau) = f \cdot Q^\downarrow(\tau) + \gamma [p_0 - p_{wf}(t)] \tag{3.23}$$

or

$$-c_t V_\phi J^{-1} \cdot q^\uparrow(t) + Q^\uparrow(\tau) = f \cdot Q^\downarrow(\tau) + c_t V_\phi [p_0 - p_{wf}(t)] \tag{3.24}$$

where $\tau$ and $\gamma$ constants are related to primary well and reservoir properties:

$$\gamma = c_t V_\phi \tag{3.25}$$

$$\tau = J^{-1} \cdot \gamma = J^{-1} \cdot c_t \cdot V_\phi \tag{3.26}$$

Generalization for arbitrary number of producers and injectors is quite straightforward and leads to (3.11).

## 4. Link between MDCV and CRM

It will be proven below that **CRM** equation (3.11) can be derived from **Multiwell Convolution Equation** (2.24) with the following assumptions on **UTR**s:

| DTR | CTR from offset producers | CTR from offset injection |
|---|---|---|
| $p_{u,nn}(t) = \frac{\tau_n}{\gamma_n} \cdot \theta(t) + \frac{t}{\gamma_n}$, $n \in \uparrow$  (4.1) | $p_{u,nm}(t) = 0$, $m \in \uparrow$   (4.2) | $p_{u,nm}(t) = \frac{f_{nm}}{\gamma_n} \cdot t$, $m \in \downarrow$   (4.3) |

The physical meaning of (4.1) is that $n$-th producer is draining a limited volume area with no wellbore storage and immediate (after the flow starts) boundary-dominated pseudo-state pressure diffusion regime.

The physical meaning of (4.2) is that $n$-th producer is not interfering with any other producer. The physical meaning of (4.3) is that $m$-th injector is sharing $f_{nm}$ portion of its total injection towards supporting a pressure in $n$-th producer and any change in injection translates to the immediate pressure response in $n$-th producer with boundary-dominated pseudo-state regime.

The above assumptions provide a clear insight (see below) into the limitations which **CRM** imposes on flow regimes and cross-well interference behaviour as against the more general **MDCV** case.

The equations (4.1)–(4.4) lead to the following pressure derivatives:

| | | |
|---|---|---|
| $\dot{p}_{u,nn}(t) = \frac{\tau_n}{\gamma_n} \cdot \delta(t) + \frac{1}{\gamma_n}$, $n \in \uparrow$ | $\dot{p}_{u,nm}(t) = 0$, $m \in \uparrow$   (4.4) | $\dot{p}_{u,nm}(t) = \frac{f_{nm}}{\gamma_n}$, $m \in \uparrow$ |

where $\delta(t)$ is Dirac delta function:

$$\delta(t) = \frac{d\theta(t)}{dt} \qquad (4.5)$$

Substituting the above UTRs to the convolution equation (2.24) leads to:

$$p_n(t) = p_{n0} - \int_0^t \left[\frac{\tau_n}{\gamma_n}\delta(t-\tau) + \frac{1}{\gamma_n}\right] \cdot q_n(\tau)d\tau - \sum_{m \in \downarrow} \int_0^t \frac{f_{nm}}{\gamma_n} \cdot q_m(\tau)d\tau \qquad (4.6)$$

or

$$p_n(t) = p_{n0} - \frac{\tau_n}{\gamma_n}\int_0^t \delta(t-\tau) \cdot q_n(\tau)d\tau - \frac{1}{\gamma_n} \cdot \int_0^t q_n(\tau)d\tau - \sum_{m \in \downarrow} \frac{f_{nm}}{\gamma_n} \int_0^t q_m(\tau)d\tau \qquad (4.7)$$

or

$$p_n(t) = p_{n0} - \frac{\tau_n}{\gamma_n} q_n(t) - \frac{1}{\gamma_n} \cdot \int_0^t q_n(\tau)d\tau - \sum_{m \in \downarrow} \frac{f_{nm}}{\gamma_n} \int_0^t q_m(\tau)d\tau \qquad (4.8)$$

or

$$\gamma_n \cdot [p_n(t) - p_{n0}] = -\tau_n q_n(t) - \int_0^t q_n(\tau)d\tau - \sum_{m \in \downarrow} f_{nm} \int_0^t q_m(\tau)d\tau \qquad (4.9)$$

where index $_m$ is running through injectors only: $m \in \downarrow$.

Let's define cumulative production as:

$$Q_n^\uparrow(t) = \int_0^t q_n(\tau)d\tau \tag{4.10}$$

so that

$$q_n(t) = \frac{dQ_n^\uparrow}{dt} \tag{4.11}$$

Let's also define injection rates with positive values:

$$q_m^\downarrow(t) = -q_m(t) \geq 0 \tag{4.12}$$

and corresponding injection cumulative:

$$Q_m^\downarrow(t) = -Q_m(t) \geq 0 \tag{4.13}$$

The deconvolution equation (4.9) becomes:

$$\gamma_n \cdot [p_n(t) - p_{n0}] = -\tau_n q_n(t) - Q_n^\uparrow(t) + \sum_{m \in \downarrow} f_{nm} Q_m^\downarrow(t) \tag{4.14}$$

or

$$\tau_n \frac{dQ_n^\uparrow}{dt} + Q_n^\uparrow(t) = \sum_{m \in \downarrow} f_{nm} Q_m^\downarrow(t) + \gamma_n \cdot [p_{n0} - p_n(t)] \tag{4.15}$$

which is the base form of **ICRM** (3.11).

It should be noted that the above exercise does not strictly mean that **CRM** will provide the same results as **MDCV** even in the case when **CRM** assumptions are fully met. The computational engines of **CRM** and **MDCV** are quite different and may arrive to different solutions on the same dataset, with **MDCV** having more flexibility to tune up for the dataset. Anyway, the fact that **CRM** equations can be derived from **Multiwell Convolution Equation** helps understand the scope of **CRM** limitations which are additional to **MDCV**.

## 5. Scope of MDCV applicability in production-injection pairing

In general, **MDCV** is applicable for all types of production: natural depletion, waterflood and gas flood. But since this paper sets focus on the link between **MDCV** and **CRM** and the latter are normally defined for production-injection pairing this paragraph will customise **MDCV** applicability to waterflood applications only.

It should be noted that all statements and estimates provided below are the results of the extensive engineering practice using the **Polygon** software [18] which features its own implementation of **MDCV** and **CRM** algorithms. Using other solutions may lead to a different experience.

1. Rates should vary both in number of events and in amplitude of variations, specifically:

   - Production rates should vary in order to pick up a correlation between production rate and bottomhole pressure in producer and assess producer's Drawdown Transient Response
   - Injection rates should vary in order to pick up a correlation between injectors and producers and assess injector-producer Cross-well Transient Response

2. Rates should be cumulatively accurate and the deviation from the true values should be within 30% range. The rate errors normally do not affect the convergence of the algorithm but may mislead on the **UTR**s, particularly on quantitative assessment of transmissibility. The systematic mistakes in rate records are the worse as **MDCV** has no means to understand that.

3. PDG readings should be accurate (the usual requirement is to provide 1 bar resolution or better and sampling rate of one hour or more frequent)

4. Transient responses should stay constant in time which happens when well-reservoir contacts (skin-factor, perforation interval, lateral leg length, hydraulic fracture length etc) do not change.

5. The pressure diffusion should honour the linear diffusion equation, which means that:

   - Modelling period is short enough to assume that reservoir saturation does not change substantially

   - In case of light oil and free gas the formation pressure does not change much during the modelling period otherwise the pseudo-variables should be used in **MDCV**

Sometimes diffusion is not stable in time and **UTR**s change their properties (as a result of workover, or saturation change or arise of unwanted communication or skin build-up). There are modifications to **MDCV** to account for the varying skin, providing that pressure log-derivative stays constant anyway, but it falls outside the scope of this paper.

In many practical cases, the production history can be split into the time segments where **UTR**s can be considered relatively constant and then apply **MDCV** procedure independently for each time segment. In author's experience the most frequent **MDCV** failures (in both synthetic and field cases) were related to the low number and poor amplitude of rate variation events.
It usually happens in three cases:

1. A properly compensated steady state production with constant bottomhole pressure and rate

2. Short production/injection history (few months only) with no major events

3. Gradually changing skin-factor which requires splitting the global time frame into segments which may not hold enough rate variations to facilitate the training

Sometimes one has to deal with very poor PDG readings, both in terms of resolution (> 1 bar) and sampling rate (one point a day or more). This requires very long historical records (few years) for meaningful **MDCV** analysis and may not capture efficiently the weak interference from the remote wells. But even in this case the results of **MDCV** analysis help analyst to assess **DTR** and interpret it as a traditional long-term drawdown test, cleaned from interference with offset wells. This actually makes **MDCV** a very efficient tool for rate allocation between wells which is important inout for future flow-simulation modelling.

The next popular issue is inaccurate daily and/or monthly rate allocation based on rare well tests. The commercial implementations of the **MDCV** algorithm are trying to correct the rates to honour the constant transient response principle and match the monthly multi-well cumulative (if the latter option selected by analyst).

The **PolyGon MDCV** proved the ability to correct the non-systematic errors with up to 30 % value and up to few hours of desynchronization. These are the average numbers over the typical one year dataset and may vary from case to case. The best rate correction is achieved by **MDCV** when apart from the daily allocated flowrate history the user additionally uploads reliable well test data as calibration points and reliable data on monthly cumulative on the tested group of wells. The nonlinear diffusion effects are effectively compensated by using the **MDCV** over the pseudo-pressure/pseudo-time and rate-dependent skin-factor.

## 6. Scope of CRM applicability in production-injection pairing

Conventional **CRM** is only modelling production-injection pairing and cannot be explicitly applied to interference between two producers or two injectors as against **MDCV** which can handle cross-well interference between all types of wells.

The **Polygon** software uses both conventional **CRM** (with producer-injector pairing) and modified version of **CRM** which includes producer-producer interference (as linear correlation, much in the same way as injector-producer pairing).

The modified version of **CRM** extends its applicability to the cases where producer stay in a strong interference with each other, but it does not change the scope of **CRM** applicability below.
The **CRM** scope of applicability includes all the **MDCV** constraints and additional limitations:

1. The dynamic drainage volume of producers should stay constant over time.
   It disqualifies the cases when production rates are varying substantially or when some producers are shut-in for prolonged periods of time. This clause is also violated when producer has monotonous increase or decrease in production as it may lead to substantial change in drainage volume.
   Unlike **CRM**, the **MDCV** comprehends arbitrary variations of dynamic drainage volumes around producers as soon as there are enough events in production history of surrounding producers to deconvolve their connectivity.

2. The production should be strictly pseudo steady-state
   It disqualifies the cases when producer is draining very large volume or when production is efficiently supported by gas, aquifer or unaccounted offset injection (like those caused by thief injection in offset injector from another pay).
   Unlike **CRM**, the **MDCV** can successfully deconvolve and model the pressure diffusion under various boundary conditions: no-flow boundaries, infinite reservoir and full/partial pressure support.
   This gives **MDCV** a clear benefit of deciphering the well performance, the cross-well connectivity and additionally – to assess the boundary conditions outside the group of tested wells and compare this with expectations from 3D model.

3. The productivity index of producers should stay constant over time
   It disqualifies the cases when transient production (which sometimes is called "fill-in period") constitutes a large portion of total production period and productivity index is constantly varying.
   It usually happens in mid to low permeability reservoirs or heavy oil reservoirs.
   Unlike **CRM**, the **MDCV** model is using the full-range unit-rate transient responses and properly accounts for the transition effects and varying productivity index.

4. Conventional **CRM** does not work for channel flow and near-boundary flow due to substantial deviation of **DTR** from the unit-slope.
   Unlike **CRM**, the **MDCV** model is using the flexible **DTR** and can tune up for all types of boundaries.

5. Conventional **CRM** shows poor predictions on horizontal producers due to substantial deviation of horizontal **DTR** from the unit-slope **DTR** assumed by **CRM**.
   Unlike **CRM**, the **MDCV** model is using the flexible **DTR** and can tune up for fractured wells, both vertical and horizontals.

6. Conventional **CRM** does not work on strongly compressible fluids, although there are some academic studies in this area.
   Unlike **CRM**, the **MDCV** model can be applied to pseudo-potentials and account for the skin-rate dependence thus covering a wide range of non-linear effects of strongly compressible fluids.

The ultimate majority of **CRM** failures in author's experience have found to be violating the above principles.

## 7. Assessing the model prediction quality

All the above limitations for **MDCV** / **CRM** applicability are not constructive and do not provide a practical guideline for the analyst to assess the model applicability in each case.

The simplest way to assess the prediction quality of training process is to perform cross-validation. The initial historical database is manually split into two sub-sets: training and validation, both having at least few noticeable events.

Once the training is completed the model prediction are tested over the validation dataset and the deviation with historical data is reported as the prediction quality (usually in terms of coefficient of determination $R^2$ and root-mean-square deviation RMSD).

The **MDCV** model has some substantial advantage over **CRM** in terms of cross-validation. Since **MDCV** is capable to pick up a transient behavior of pressure diffusion it can be validated against the short-term pressure build-ups, fall-offs or drawdowns, while **CRM** should be trained and validated on boundary-dominated drawdown production regimes only. This may not help if the wells were not producing much time in transition or when build-up behavior does not match the drawdown behavior (like in case of dynamic fractures).

The selection of training vs validation dataset is not always obvious, especially if the history is quite short and uneventful and one cannot eye-pick two subsets both having substantial events.
The higher the validation/training dataset volume ratio the higher the trust to the predictive quality of the model. Although the sheer volume is not the only qualifier and the amplitude of the pressure-rate variations have also a great effect on the predictive quality of the model. If the choice of the split is obvious and one arrives to a strong cross-validated solution on the eventful validation dataset then the choice was ok and there is no need in further endeavours. Otherwise, one need to be sure that the training-validation split is fair.

The most practical solution which author was using for a long time in this case is to:

1. Build a simple numerical pressure diffusion PEBI-grid model around **MDCV/CRM** wells
    a. Define model boundaries
    b. Place the well trajectories
    c. Populate representative fluid, well and reservoir properties
    d. Distribute the properties assuming connectivity between all wells
    e. Upload the true production/injection histories for all wells from the field data records
2. Generate the synthetic bottom-hole pressure history which will account for the typical pressure-related non-linearities and saturation dynamics (assuming all wells are connected)
3. Upload the flowrate-pressure dataset to **MDCV/CRM** modelling facility
4. Split the database into training and validation datasets and select the training-engine properties
5. Train the model with training dataset
6. Cross-validate its rate or pressure predictions with validation dataset
7. If not validated (with criterion being set by user depending on the purposes of **MDCV/CRM**) – then go back to step 4 and repeat the 5+6 with the new split
8. If few iterations do not provide a fair cross-validation, then report "no-validation conditions" and either discard the case or perform **MDCV/CRM** training without cross-validation (which means low trust on the end result)
9. Otherwise, use the last training-validation dataset split and perform **MDCV/CRM** training with cross-validation using the actual flowrate-pressure dataset from the field records
10. If validation is not achieved, then report the **MDCV/CRM** failure

The above workflow is using synthetic modelling only to assess the fair training/validation dataset split or raise concern on applicability of **MDCV/CRM** on this dataset. The simulated pressure itself has no meaning in this exercise. If the data processing facility is up to the job (the author is using **Polygon** pressure/rate modelling facility [18] and includes PEBI-grid multi-well simulations, **MDCV** and **CRM**) then all the above actions take less than an hour which is acceptable in most practical applications and does not add much to the total **MDCV/CRM** turn-around-time.

## 8. Conclusions

Both **MDCV** and **CRM** are an effective production analysis tools which can realistically model how injection impacts the liquid production and can be a part of production optimisation planning.

They also provide calibration data points for 3D modelling, which means that the ultimate full-field flow simulation 3D model should match the **UTR**s from **MDCV** or cross-well interference coefficients and productivity index from **CRM**.

The main limitation of **MDCV** and **CRM** is that flowrates must vary (both in the number of events and the amplitude of variation) in order to train the model parameters as uniquely as possible.

There is no simple efficient metric to assess if a given flowrate-pressure history can be decoded by **MDCV/CRM** and it requires of **MDCV/CRM** to be actually executed and check if the post-factum cross-validation is successful.

The only exception is when rates do not vary at all which is a straightforward indication that both **MDCV** and **CRM** will fail.

If splitting the history into training/validation is not obvious then author's recommendation is to perform simple synthetic PEBI-grid exercises with true rate variation history and estimate the fair training/validation split and if successful, then proceed into the **MDCV/CRM** with actual field bottomhole pressure data.

There are specific assumptions of **CRM** which limit its use comparing to **MDCV**:

- The geological or dynamic boundary of the tested group of wells should be "no-flow" resulting in pseudo-steady state flow regime for all producers if all injection is shut-in.

- The dynamic drainage areas of producers should stay constant during the time period at which the model is trained and used for production forecast. This particularly dictates that production rates should not have any ascending/descending trends.

- The permeability should be high enough to provide fast transition to the stabilised flow (with constant productivity index) so that most of the time the wells are producing at stabilised regime and not in transition.

- The conventional **CRM** engines do not perform rate corrections which means that the injection rate histories should be accurate to achieve a cross-validated training of model parameters.

- **CRM** modelling should be avoided in long horizontals (usually > 500 m) and heavily fractured vertical wells.

## Acknowledgments

Author would like to thank prof. Marat Timerbaev for his help in reading this paper and making a fair point on the expanding the class of the boundary value problems which allow convolution principle.